\documentclass{amsart}%
\usepackage{graphicx}
\usepackage{amscd}
\usepackage{amsmath}
\usepackage{amsfonts}
\usepackage{amssymb}%
\newtheorem{theorem}{Theorem}
\theoremstyle{plain}

\newtheorem{corollary}{Corollary}

\newtheorem{proposition}{Proposition}
\newtheorem{remark}{Remark}

\numberwithin{equation}{section}

\begin{document}
\title[Multilinear extensions]{On the multilinear extensions of the concept of absolutely summing operators}
\author{Daniel M. Pellegrino}
\address[Daniel M. Pellegrino]{Depto de Matem\'{a}tica e Estat\'{\i}stica- Caixa Postal
10044- UFCG- Campina Grande-PB-Brazil }
\email{dmp@dme.ufcg.edu.br}
\thanks{This work is partially supported by Instituto do Mil\^{e}nio, IMPA}

\begin{abstract}
In this paper we investigate the several different extensions of the concept
of absolutely summing operators and their connections.

\end{abstract}
\maketitle

\section{Introduction and notation}

The core of the theory of absolutely summing operators lie in the ideas of A.
Grothendieck in the 1950s. Further work (after a decade) of A. Pietsch
\cite{P1967} and Lindenstrauss and Pe\l \ czy\'{n}ski \cite{Lindenstrauss}
clarified Grothendieck$^{\prime}$s insights and nowadays the ideal of
absolutely summing operators is a central topic of investigation. For details
on absolutely summing operators we refer to the book by Diestel-Jarchow-Tonge
\cite{Diestel}.

A natural question is how to extend the concept of absolutely summing
operators to multilinear mappings and polynomials. A first light in this
direction is the work by Alencar-Matos \cite{AlencarMatos}, where several
classes of multilinear mappings between Banach spaces were investigated. Since
then, just concerning to the idea of lifting the ideal of absolutely summing
operators to polynomials and multilinear mappings, there are several works in
different directions (we mention Bombal \textit{et al} \cite{Perez}, Dimant
\cite{Dimant}, Matos \cite{Anais},\cite{Collect},\cite{Nach}). However, there
seems to be no effort in the direction of comparing these different classes.
The aim of this paper is to investigate these classes and their connections.

Throughout this paper $E,E_{1},...,E_{n},G_{1},...,G_{n,}F,F_{0}$ will be
Banach spaces. Given a natural number $n\geq2,$ the Banach space of all
continuous $n$-linear mappings from $E_{1}\times...\times E_{n}$ into $F$
endowed with the $\sup$ norm will be denoted by $\mathcal{L}(E_{1},...,E_{n}%
$;$F)$ and the space of all continuous $n$-homogeneous polynomials $P$ from
$E$ into $F$ with the $\sup$ norm is represented by $\mathcal{P}(^{n}E;F).$ If
$T$ is a multilinear mapping and $P$ is the polynomial generated by $T$, we
write $P=\overset{\wedge}{T}.$ Conversely, for the (unique) symmetric
$n$-linear mapping associated to an $n$-homogeneous polynomial $P$ we use the
symbol $\overset{\vee}{P}.$ For $i=1,...,n,$ we denote by $\Psi_{i}%
^{(n)}:\mathcal{L}(E_{1},...,E_{n};F)\rightarrow\mathcal{L}(E_{i}%
;\mathcal{L}(E_{1},\overset{[i]}{...},E_{n};F))$ the canonical isometric
isomorphism
\[
\Psi_{i}^{(n)}(T)(x_{i})(x_{1}\overset{[i]}{...}x_{n})=T(x_{1},...,x_{n}),
\]
where $\overset{[i]}{...}$ means that the $i$-th coordinate is not involved.

For $p\in]0,\infty\lbrack,$ the linear space of all sequences $(x_{j}%
)_{j=1}^{\infty}$ in $E$ such that $\Vert(x_{j})_{j=1}^{\infty}\Vert_{p}%
=(\sum_{j=1}^{\infty}\Vert x_{j}\Vert^{p})^{\frac{1}{p}}<\infty$ is denoted by
$l_{p}(E).$ We represent by $l_{p}^{w}(E)$ the linear space of the sequences
$(x_{j})_{j=1}^{\infty}$ in $E$ such that $(\varphi(x_{j})))_{j=1}^{\infty}\in
l_{p}$ for every continuous linear functional $\varphi:E\rightarrow
\mathbb{K},$ and define $\Vert.\Vert_{w,p}$ in $l_{p}^{w}(E)$ by $\Vert
(x_{j})_{j=1}^{\infty}\Vert_{w,p}=\sup_{\varphi\in B_{E}%
\acute{}%
}\Vert(\varphi(x_{j}))_{j=1}^{\infty}\Vert_{p}.$ If $p=\infty$ we are
restricted to the case of bounded sequences and in $l_{\infty}(E)$ we use the
$\sup$ norm. One can verify that $\Vert.\Vert_{p}$ $(\Vert.\Vert_{w,p})$ is a
$p$-norm in $l_{p}(E)($ $l_{p}^{w}(E))$ for $p<1$ and a norm in $l_{p}(E)($
$l_{p}^{w}(E))$ for $p\geq1.$

We begin by presenting the several classes of multilinear mappings related to
the concept of absolutely summing operators:

\begin{itemize}
\item $T\in\mathcal{L}(E_{1},...E_{n};F)$ is said to be $p$\emph{-dominated}
if there exist $C\geq0$ and regular probability measures $\mu_{j}$ on the
Borel $\sigma$-algebras $\mathcal{B}(B_{E_{j}^{^{\prime}}})$ of $B_{E_{j}%
^{^{\prime}}}$ endowed with the weak star topologies $\sigma(E_{j}^{\prime
},E_{j}),$ $j=1,...,n,$ such that
\[
\left\Vert T\left(  x_{1},...,x_{n}\right)  \right\Vert \leq C\prod
\limits_{j=1}^{n}\left[  \int_{B_{E_{j}^{\prime}}}\left\vert \varphi\left(
x\right)  \right\vert ^{p}d\mu_{j}\left(  \varphi\right)  \right]  ^{\frac
{1}{p}}%
\]
for every $x_{j}\in E_{j}$ and $j=1,...,n$. It is well known that
$T\in\mathcal{L}(E_{1},...E_{n};F)$ is $p$-dominated if and only if there
exist Banach spaces $G_{1},...,G_{n},$ absolutely $p$-summing linear operators
$u_{j}\in\mathcal{L}(E_{j};G_{j})$ and a continuous $n$-linear mapping
$R\in\mathcal{L}(G_{1},...,G_{n};F)$ so that $T=R(u_{1},...u_{n})$. Notation:
$T\in\mathcal{L}_{d,p}(E_{1},...E_{n};F).$

\item $T\in$ $\mathcal{L}(E_{1},...,E_{n};F)$ is said to be \emph{of
absolutely }$p$\emph{-summing type }$(T\in$ $[\Pi_{as(p)}](E_{1}%
,...,E_{n};F))$ if $\Psi_{j}^{(n)}(T)$ is absolutely $p$-summing for every $j$
$\in\{1,...,n\}.$

\item $T\in\mathcal{L}(E_{1},...E_{n};F)$ is\emph{ }$p$\emph{-semi-integral}
if there exist $C\geq0$ and a regular probability measure $\mu$ on the Borel
$\sigma-$algebra $\mathcal{B}(B_{E_{1}^{^{\prime}}}\times...\times$
$B_{E_{n}^{^{\prime}}})$ of $B_{E_{1}^{^{\prime}}}\times...\times$
$B_{E_{n}^{^{\prime}}}$ endowed with the weak star topologies $\sigma
(E_{l}^{\prime},E_{l}),$ $l=1,...,n,$ such that
\[
\parallel T(x_{1},...,x_{n})\parallel\leq C\left(  \int_{B_{E_{1}^{\prime}%
}\times...\times B_{E_{n}^{\prime}}}\mid\varphi_{1}(x_{1})...\varphi_{n}%
(x_{n})\mid^{p}d\mu(\varphi_{1},...,\varphi_{n})\right)  ^{\frac{1}{p}}%
\]
for every $x_{j}\in E_{j}$ and $j=1,...,n$. Notation: $T\in\mathcal{L}%
_{si,p}(E_{1},...E_{n};F).$ The infimum of the $C$ defines a norm $\left\Vert
.\right\Vert _{si,p}$ for the space of $p$-semi integral mappings.

\item $T:E_{1}\times...\times E_{n}\rightarrow F$ is \emph{fully (or multiple)
}$p$\emph{-summing} if there exists $C>0$ such that
\[
\left(  \sum_{j_{1},...,j_{n}=1}^{\infty}\Vert T(x_{j_{1}}^{(1)},...,x_{j_{n}%
}^{(n)})\Vert^{p}\right)  ^{\frac{1}{p}}\leq C\prod\limits_{k=1}^{n}%
\Vert(x_{j}^{(k)})_{j=1}^{\infty}\Vert_{w,p}\text{ }\forall(x_{j}^{(k)}%
)_{j=1}^{\infty}\in l_{p}^{w}(E_{k}),k=1,...,n.
\]
The space of all fully $p$-summing $n$-linear mappings from $E_{1}%
\times...\times E_{n}$ into $F$ will be denoted by $\mathcal{L}_{fas,p}%
(E_{1},...,E_{n};F),$ and the infimum of the $C$ for which the inequality
always holds defines a norm $\Vert.\Vert_{fas,p}$ for $\mathcal{L}%
_{fas,p}(E_{1},...,E_{n};F)$.
\end{itemize}

\begin{itemize}
\item $T\in\mathcal{L}(E_{1},...E_{n};F)$ is strongly $p$-summing if there
exists $C\geq0$ and a regular probability measure $\mu$ on the Borel $\sigma
-$algebra $\mathcal{B}(B_{\mathcal{L}(E_{1},...,E_{n};\mathbb{K})})$ of
$B_{\mathcal{L}(E_{1},...,E_{n};\mathbb{K})}$ with the weak star topology such
that
\[
\parallel T(x_{1},...,x_{n})\parallel\leq C\left(  \int_{B_{\mathcal{L}%
(E_{1},...,E_{n};\mathbb{K})}}\mid\phi(x_{1},...,x_{n})\mid^{p}d\mu
(\phi)\right)  ^{\frac{1}{p}}%
\]
for every $x_{j}\in E_{j}$ and $j=1,...,n$. Notation: $T\in\mathcal{L}%
_{sas,p}(E_{1},...E_{n};F):$

\item $T\in\mathcal{L}(E_{1},...E_{n};F)$ is \emph{absolutely }$(p;q_{1}%
,...,q_{n})$\emph{-summing} \emph{(or }$(p;q_{1},...,q_{n})$\emph{-summing) at
the point }$(a_{1},...,a_{n})$ in $E_{1}\times...\times E_{n}$ if
\[
(T(a_{1}+x_{j}^{(1)},...,a_{n}+x_{j}^{(n)})-T(a_{1},...,a_{n}))_{j=1}^{\infty
}\in l_{p}(F)
\]
for every $(x_{j}^{(s)})_{j=1}^{\infty}\in l_{q_{s}}^{w}(E)$, $s=1,...,n.$ In
the case that $T$ is $(p;q_{1},...,q_{n})$-summing at every $(a_{1}%
,...,a_{n})\in E_{1}\times...\times E_{n}$ we say that $T$ is $(p;q_{1}%
,...,q_{n})$-summing everywhere. Notation: $T\in\mathcal{L}_{as(p,q_{1}%
,...,q_{n})}^{ev}(E_{1},...E_{n};F).$ If $T$ is $(p;q_{1},...,q_{n})$-summing
at $(0,...,0)\in E_{1}\times...\times E_{n}$ we say that $T$ is $(p;q_{1}%
,...,q_{n})$\emph{-summing }and we write $T\in\mathcal{L}_{as(p,q_{1}%
,...,q_{n})}(E_{1},...E_{n};F)$. When $p=q_{1},...,q_{n}$ we write
$\mathcal{L}_{as,p}^{ev}(E_{1},...E_{n};F)$ and/or $\mathcal{L}_{as,p}%
(E_{1},...E_{n};F).$ It is well known that $\mathcal{L}_{d,p}(E_{1}%
,...E_{n};F)=\mathcal{L}_{as(\frac{p}{n},p,...,p)}(E_{1},...E_{n};F)$.
\end{itemize}

Except perhaps for the concept of $p$-semi integral mappings, all of the above
concepts are well known and individually investigated. The $p$-semi integral
mappings were motivated by the work of Alencar-Matos \cite{AlencarMatos} and
introduced in \cite{Pellegrino}. The dominated mappings were first explored by
Schneider \cite{Schneider} and Matos \cite{Anais} and more recently in
\cite{Botelhoirish},\cite{BP},\cite{Melendez} and \cite{studia}. Multilinear
mappings of absolutely summing type are motivated by abstract methods of
creating ideals and are explored in \cite{Pellegrino}. The ideal of fully
summing multilinear mappings was introduced by Matos \cite{Collect} and
investigated by M. L. Souza \cite{Souza} in her doctoral thesis under his
supervision. It was also independently introduced by Bombal \textit{et al
}(with a different name \textquotedblleft multilple summing\textquotedblright)
and explored in \cite{Perez}. The ideal of strongly summing multilinear
mappings was introduced by V. Dimant \cite{Dimant} and the absolutely summing
multilinear mappings appears in \cite{AlencarMatos},\cite{Anais} and have been
vastly studied (we mention \cite{Nach},\cite{studia},\cite{Pellegrino}%
,\cite{cotype},\cite{PG} for example). In the next two sections we investigate
the $p$-semi integral and absolutely summing mappings. In Section 4 we study
the connections between the classes previously introduced and in the last
section we define a new related class and sketch their main properties.

\section{$p$-semi integral mappings}

We begin with a characterization of $p$-semi-integral mappings, that will be
useful in section 4:

\begin{theorem}
\label{02}$T\in\mathcal{L}_{si,p}(E_{1},...E_{n};F)$ if and only if there
exists $C\geq0$ such that
\begin{equation}
\left(  \sum\limits_{j=1}^{m}\parallel T(x_{1,j},...,x_{n,j})\parallel
^{p}\right)  ^{1/p}\leq C\left(  \underset{l=1,...,n}{\underset{\varphi_{l}\in
B_{E_{l}^{\prime}}}{\sup}}\sum\limits_{j=1}^{m}\mid\varphi_{1}(x_{1,j}%
)...\varphi_{n}(x_{n,j})\mid^{p}\right)  ^{1/p} \label{day}%
\end{equation}
for every natural $m,$ $x_{l,j}\in E_{l}$ with $l=1,...,n$ and $j=1,...,m.$ We
also have that the infimum of the $C$ is $\parallel T\parallel_{si,p}.$
\end{theorem}

Proof. If $T$ is $p$-semi integral, it is not hard to obtain (\ref{day}).

Conversely, suppose that (\ref{day}) holds. The proof follows the idea of the
case $p=1$ in \cite{AlencarMatos}. Define

\begin{itemize}
\item $\Gamma_{1}=\{f\in C(B_{E_{1}^{\prime}}\times...\times B_{E_{n}^{\prime
}});f<C^{-p}\}.$

\item $\Gamma_{2}=co\{f\in C(B_{E_{1}^{\prime}}\times...\times B_{E_{n}%
^{\prime}});$ there are $x_{l}\in E_{l},$ $l=1,...,n,$ so that $\parallel
T(x_{1},...,x_{n})\parallel=1$ and $f(\varphi_{1},...\varphi_{n})=\mid
\varphi_{1}(x_{1})...\varphi(x_{n})\mid^{p}\}.$
\end{itemize}

where $co\{.\}$ denotes the convex hull. Let us show that $\Gamma_{1}%
\cap\Gamma_{2}=\phi.$

If $h\in\Gamma_{2},$ then $h=\sum\limits_{j=1}^{m}\alpha_{j}f_{j},$
$\alpha_{j}>0,\sum\limits_{j=1}^{m}\alpha_{j}=1$ and
\[
f_{j}(\varphi_{1},...\varphi_{n})=\mid\varphi_{1}(x_{1,j})...\varphi
_{n}(x_{n,j})\mid^{p}\text{ \ }\forall\varphi_{l}\in B_{E_{l^{^{\prime}}}}.
\]
By hypothesis we have
\[
\left\Vert h\right\Vert =\left(  \underset{l=1,...,n}{\underset{\varphi_{l}\in
B_{E_{l}^{\prime}}}{\sup}}\sum\limits_{j=1}^{m}\mid\varphi_{1}(\alpha
_{j}^{1/p}x_{1,j})...\varphi_{n}(x_{n,j})\mid^{p}\right)  \geq C^{-p}%
\sum_{j=1}^{m}(\alpha_{j}^{1/p})^{p}\parallel T(x_{1,j},...,x_{n,j}%
)\parallel^{p}=C^{-p}.
\]
Hence $h\notin\Gamma_{1}.$ By Hahn-Banach Theorem there exist $\lambda>0$ and
\[
\psi\in C(B_{E_{1}^{\prime}}\times...\times B_{E_{n}^{\prime}})^{\prime}%
\]
so that $\parallel\psi\parallel=1$ and
\[
\psi(f)\leq\lambda\leq\psi(g)\text{ }\forall f\in\Gamma_{1},g\in\Gamma_{2}.
\]
Since each $\ f<0$ belongs to $\Gamma_{1},$ we have $\psi(mf)\leq\lambda$ for
every natural $m.$ Thus $\psi(f)\leq0$ and $\psi$ is a positive functional and
thus there exists a regular probability measure $\mu,$ defined on the Borel
sets of $B_{E_{1}^{\prime}}\times...\times B_{E_{n}^{\prime}}$ (with the weak
star topology)$,$ so that
\[
\psi(f)=\int_{B_{E_{1}^{\prime}}\times...\times B_{E_{n}^{\prime}}}fd\mu.
\]
Defining $f_{m}$ by $f_{m}=C^{-p}-\frac{1}{m},$ we have $f_{m}\in\Gamma_{1}$
for every natural $m$. Thus
\[
\int_{B_{E_{1}^{\prime}}\times...\times B_{E_{n}^{\prime}}}f_{m}d\mu
=C^{-p}-\frac{1}{m}\leq\lambda\text{ for every }m,
\]
and hence $\lambda\geq C^{-p}.$

Therefore, if $\parallel T(x_{1},...,x_{n})\parallel=1,$ defining
$f(\varphi_{1},...\varphi_{n}):=\mid\varphi_{1}(x_{1})...\varphi(x_{n}%
)\mid^{p},$ we have $f\in\Gamma_{2}$ and%
\begin{equation}
\int_{B_{E_{1}^{\prime}}\times...\times B_{E_{n}^{\prime}}}fd\mu=\psi(f)\geq
C^{-p}=C^{-p}\parallel T(x_{1},...,x_{n})\parallel\label{az6}%
\end{equation}
i.e.,
\[
C^{p}\int_{B_{E_{1}^{\prime}}\times...\times B_{E_{n}^{\prime}}}\mid
\varphi_{1}(x_{1})...\varphi(x_{n})\mid^{p}d\mu\geq\parallel T(x_{1}%
,...,x_{n})\parallel,
\]
and since $\parallel T(x_{1},...,x_{n})\parallel=1$, we obtain
\[
\parallel T(x_{1},...,x_{n})\parallel\leq C\left(  \int_{B_{E_{1}^{\prime}%
}\times...\times B_{E_{n}^{\prime}}}\mid\varphi_{1}(x_{1})...\varphi
(x_{n})\mid^{p}d\mu\right)  ^{\frac{1}{p}}.
\]
If $\parallel T(x_{1},...,x_{n})\parallel\neq0,$ it suffices to replace
$x_{1}$ by $x_{1}\parallel T(x_{1},...,x_{n})\parallel^{-1}$in (\ref{az6}),
and we obtain
\[
\parallel T(x_{1},...,x_{n})\parallel\leq C\left(  \int_{B_{E_{1}^{\prime}%
}\times...\times B_{E_{n}^{\prime}}}\mid\varphi_{1}(x_{1})...\varphi
(x_{n})\mid^{p}d\mu\right)  ^{\frac{1}{p}}.
\]

We list some interesting properties of $p$-semi integral mappings, whose proof
are standard and we omit:

\begin{proposition}
\label{demmm1}

(i) If $p\leq q,$ then $\mathcal{L}_{si,p}(E_{1},...E_{n};F)\subset
\mathcal{L}_{si,q}(E_{1},...E_{n};F)$.

(ii) If $T\in\mathcal{L}(E_{1},...E_{n};F)$ is $p$-semi integral and
$i:F\rightarrow F_{0}$ is an isometric embedding, then $i\circ T$ is $p$-semi
integral and $\parallel i\circ T\parallel_{si,p}=\parallel T\parallel_{si,p}.$

(iii) If $T\in\mathcal{L}_{si,p}(E_{1},...,E_{n};F),$ then $\Psi_{i}%
^{(n)}(T)\in\mathcal{L}_{as,p}(E_{i;}\mathcal{L}(E_{1},\overset{[i]}%
{...},E_{n};F))$ and $\Psi_{i}^{(n)}(T)(x)$ is $p$-semi-integral for every $x$
in $E_{i}.$

(iv) If $\mathcal{L}(E_{1},...E_{n};F)=\mathcal{L}_{si,p}(E_{1},...E_{n};F),$
then $\mathcal{L}(E_{j_{1}},...E_{j_{n}};F)=\mathcal{L}_{si,p}(E_{j_{1}%
},...E_{j_{n}};F)$, for every $j_{1},...,j_{k}$ in $\{1,...,n\}$ with
$j_{r}\neq j_{s}$ for $r\neq s.$
\end{proposition}

\section{\bigskip Absolutely summing mappings}

If we look for coincidence results, i.e., situations in which one of the
aforementioned classes coincides with the whole space of continuous
multilinear mappings, it is interesting to work with the class of absolutely
summing mappings. In \cite{Botelhoirish}, it is shown that every continuous
bilinear form defined in $\mathcal{L}_{\infty}$-spaces is absolutely
$(1;2,2)$-summing $(2$-dominated$).$ In the same paper it is also proved that
we can not expect another similar coincidence theorem for $p$-dominated
$n$-linear mappings, with $n>2$. Recently, using a generalized Grothendieck's
inequality, P\'{e}rez-Garc\'{\i}a \cite{Perez} obtained the following result
of coincidence:

\begin{theorem}
(P\'{e}rez-Garc\'{\i}a \cite{Perez})\label{tp} If $E_{1},...,E_{n}$ are
$\mathcal{L}_{\infty}$-spaces, then every continuous $n$-linear ($n\geq2$)
mapping $T:E_{1}\times...\times E_{n}\rightarrow\mathbb{K}$ is absolutely
$(1;2...,2)$-summing.
\end{theorem}

In this section we present new coincidence situations for absolutely summing
multilinear mappings. The next theorem generalizes a result theorem due to
C.A. Soares \cite{CA}:

\begin{theorem}
\label{t1}Let $A\in\mathcal{L}(E_{1},...,E_{n};F)$ and suppose that there
exists $K>0$ so that for any $x_{1}\in E_{1},....,x_{r}\in E_{r},$ the
$s$-linear ($s=n-r$) mapping $A_{x_{1}....x_{r}}(x_{r+1},...,x_{n}%
)=A(x_{1},...,x_{n})$ is absolutely $(1;q_{1},...,q_{s})$-summing and
$\left\Vert A_{x_{1}....x_{r}}\right\Vert _{as(1;q_{1},...,q_{s})}\leq
K\left\Vert A\right\Vert \left\Vert x_{1}\right\Vert ...\left\Vert
x_{r}\right\Vert $. Then $A$ is absolutely $(1;1,...,1,q_{1},...,q_{s})$-summing.
\end{theorem}

Proof. For $x_{1}^{(1)},...,x_{1}^{(m)}\in E_{1},....,x_{n}^{(1)}%
,...,x_{n}^{(m)}\in E_{n}$, let us consider $\varphi_{j}\in B_{F^{\prime}}$
such that
\[
\left\Vert A(x_{1}^{(j)},...,x_{n}^{(j)})\right\Vert =\varphi_{j}%
(A(x_{1}^{(j)},...,x_{n}^{(j)}))
\]
$\text{ for every }j=1,...,m.$ Thus, defining by $r_{j}(t)$ the Rademacher
functions on $[0,1]$ and denoting by $\lambda$ the Lebesgue measure in
$I=[0,1]^{r},$ we have
\begin{align*}
&  \int\nolimits_{I}\sum\limits_{j=1}^{m}\left(  \prod_{l=1}^{r}r_{j}%
(t_{l})\right)  \varphi_{j}A(\sum\limits_{j_{1}=1}^{m}r_{j_{1}}(t_{1}%
)x_{1}^{(j_{1})},...,\sum\limits_{j_{r}=1}^{m}r_{j_{r}}(t_{r})x_{r}^{(j_{r}%
)},x_{r+1}^{(j)},...,x_{n}^{(j)})d\lambda\\
&  =\sum\limits_{j,j_{1},...j_{r}=1}^{m}\varphi_{j}A(x_{1}^{(j_{1})}%
,...,x_{r}^{(j_{r})},x_{r+1}^{(j)},...,x_{n}^{(j)})\int\limits_{0}^{1}%
r_{j}(t_{1})r_{j_{1}}(t_{1})dt_{1}...\int\limits_{0}^{1}r_{j}(t_{r})r_{j_{r}%
}(t_{r})dt_{r}\\
&  =\sum\limits_{j=1}^{m}\sum\limits_{j_{1}=1}^{m}...\sum\limits_{j_{r}=1}%
^{m}\varphi_{j}A(x_{1}^{(j_{1})},...,x_{r}^{(j_{r})},x_{r+1}^{(j)}%
,...,x_{n}^{(j)})\delta_{jj_{1}}...\delta_{jj_{r}}\\
&  =\sum\limits_{j=1}^{m}\varphi_{j}A(x_{1}^{(j)},...,x_{n}^{(j)}%
)=\sum\limits_{j=1}^{m}\left\Vert A(x_{1}^{(j)},...,x_{n}^{(j)})\right\Vert
=(\ast).
\end{align*}
So, for each $l=1,...,r,$ assuming $z_{l}=\sum\limits_{j=1}^{m}r_{j}%
(t_{l})x_{l}^{(j)}$ we obtain
\begin{align*}
(\ast)  &  =\int\nolimits_{I}\sum\limits_{j=1}^{m}\left(  \prod_{l=1}^{r}%
r_{j}(t_{l})\right)  \varphi_{j}A(\sum\limits_{j_{1}=1}^{m}r_{j_{1}}%
(t_{1})x_{1}^{(j_{1})},...,\sum\limits_{j_{r}=1}^{m}r_{j_{r}}(t_{r}%
)x_{r}^{(j_{r})},x_{r+1}^{(j)},...,x_{n}^{(j)})d\lambda\\
&  \leq\int\nolimits_{I}\left\vert \sum\limits_{j=1}^{m}\left(  \prod
_{l=1}^{r}r_{j}(t_{l})\right)  \varphi_{j}A(\sum\limits_{j_{1}=1}^{m}r_{j_{1}%
}(t_{1})x_{1}^{(j_{1})},...,\sum\limits_{j_{r}=1}^{m}r_{j_{r}}(t_{r}%
)x_{r}^{(j_{r})},x_{r+1}^{(j)},...,x_{n}^{(j)})\right\vert d\lambda\\
&  \leq\int\nolimits_{I}\sum\limits_{j=1}^{m}\left\Vert A(\sum\limits_{j_{1}%
=1}^{m}r_{j_{1}}(t_{1})x_{1}^{(j_{1})},...,\sum\limits_{j_{r}=1}^{m}r_{j_{r}%
}(t_{r})x_{r}^{(j_{r})},x_{r+1}^{(j)},...,x_{n}^{(j)})\right\Vert d\lambda\\
&  \leq\sup_{t_{l}\in\lbrack0,1],l=1,...,r}\sum\limits_{j=1}^{m}\left\Vert
A(\sum\limits_{j_{1}=1}^{m}r_{j_{1}}(t_{1})x_{1}^{(j_{1})},...,\sum
\limits_{j_{r}=1}^{m}r_{j_{r}}(t_{r})x_{r}^{(j_{r})},x_{r+1}^{(j)}%
,...,x_{n}^{(j)})\right\Vert \\
&  \leq\sup_{t_{l}\in\lbrack0,1],l=1,...,r}\left\Vert A_{z_{1}...z_{r}%
}\right\Vert _{as(1;q_{1},...,q_{s})}\left\Vert (x_{r+1}^{(j)})_{j=1}%
^{m}\right\Vert _{w,q_{1}}...\left\Vert (x_{n}^{(j)})_{j=1}^{m}\right\Vert
_{w,q_{s}}\\
&  \leq\sup_{t_{l}\in\lbrack0,1],l=1,...,r}K\left\Vert A\right\Vert \left\Vert
z_{1}\right\Vert ...\left\Vert z_{r}\right\Vert \left\Vert (x_{r+1}%
^{(j)})_{j=1}^{m}\right\Vert _{w,q_{1}}...\left\Vert (x_{n}^{(j)})_{j=1}%
^{m}\right\Vert _{w,q_{s}}\\
&  \leq K\left\Vert A\right\Vert \left(  \prod_{l=1}^{r}\left\Vert
(x_{l}^{(j)})_{j=1}^{m}\right\Vert _{w,1}\right)  \left(  \prod_{l=1}%
^{s}\left\Vert (x_{l}^{(j)})_{j=1}^{m}\right\Vert _{w,q_{l}}\right)  .\text{ }%
\end{align*}
We have the following straightforward consequence:

\begin{corollary}
\label{ssss}If
\[
\mathcal{L}(E_{1},...,E_{m};F)=\mathcal{L}_{as(1;q_{1},...,q_{m})}%
(E_{1},...,E_{m};F)
\]
then, for any Banach spaces $E_{m+1},...,E_{n},$ we have
\[
\mathcal{L}(E_{1},...,E_{n};F)=\mathcal{L}_{as(1;q_{1},...,q_{m}%
,1,...,1)}(E_{1},...,E_{n};F).
\]

\end{corollary}

Another outcome of Theorems \ref{tp} and \ref{t1} are the following
corollaries, whose proofs are simple and we omit:

\begin{corollary}
\label{c1}If $E_{1}$,... $E_{s}$ are $\mathcal{L}_{\infty}$-spaces then, for
any choice of Banach spaces $E_{s+1},...,E_{n},$ we have
\[
\mathcal{L}(E_{1},...,E_{n};\mathbb{K})=\mathcal{L}_{as(1;q_{1},...,q_{n}%
)}(E_{1},...,E_{n};\mathbb{K}),
\]
where $q_{1}=...=q_{s}=2$ e $q_{s+1}=....=q_{n}=1.$
\end{corollary}

\begin{corollary}
\label{ff}If $\cot F=q<\infty$ \ and
\[
\mathcal{L}(E_{1},...,E_{s};\mathbb{K})=\mathcal{L}_{as(1;q_{1},....,q_{s}%
)}(E_{1},...,E_{s};\mathbb{K}),
\]
then, for any choice of Banach spaces $E_{s+1},...,E_{n},$ we have
\[
\mathcal{L}(E_{1},...,E_{n};F)=\mathcal{L}_{as(q;q_{1},....,q_{s}%
,1,....,1)}(E_{1},...,E_{n};F),
\]

\end{corollary}

\begin{corollary}
\label{g}If $\cot F=q<\infty$ and $E_{1}$,..., $E_{s}$ are $\mathcal{L}%
_{\infty}$-spaces, then, regardless of the Banach spaces $E_{s+1},...,E_{n},$
we have
\[
\mathcal{L}(E_{1},...,E_{n};F)=\mathcal{L}_{as(q;q_{1},...,q_{n})}%
(E_{1},...,E_{n};F),
\]
where $q_{1}=...=q_{s}=2$ and $q_{s+1}=...=q_{n}=1.$
\end{corollary}

It is obvious that Corollary \ref{c1} is still true if we replace $\mathbb{K}$
by any finite dimensional Banach space. A natural question is whether
Corollary \ref{c1} can be improved for some infinite dimensional Banach space
in the place of $\mathbb{K}.$ Precisely, the question is:

\begin{itemize}
\item If $E_{1},...,E_{k}$ are infinite dimensional $\mathcal{L}_{\infty}%
$-spaces, is there some infinite dimensional Banach space $F$ such that
\[
\mathcal{L}(E_{1},...,E_{k},...,E_{n};F)=\mathcal{L}_{as(1;q_{1},....,q_{n}%
)}(E_{1},...,E_{k},...,E_{n};F),
\]
where $q_{1}=...=q_{k}=2$ and $q_{k+1}=....=q_{n}=1,$ regardless of the Banach
spaces $E_{k+1},...,E_{n}$?
\end{itemize}

The answer to this question, surprisingly, is no. The proof follows directly
from \cite[Theorem 8]{studia}.

\begin{proposition}
Suppose that $E_{1},...,E_{k}$ are infinite dimensional $\mathcal{L}_{\infty}%
$-spaces. If $q_{1}=...=q_{k}=2$, $q_{k+1}=....=q_{n}=1$ and
\[
\mathcal{L}(E_{1},...,E_{k},...,E_{n};F)=\mathcal{L}_{as(1;q_{1},....,q_{n}%
)}(E_{1},...,E_{k},...,E_{n};F),
\]
regardless of the Banach spaces $E_{k+1},...,E_{n},$ then $\dim F<\infty.$
\end{proposition}

Proof. By a standard localization argument, it suffices to prove that if $\dim
F=\infty,$ then
\[
\mathcal{L}(^{n}c_{0};F)\neq\mathcal{L}_{as(1;q_{1},....,q_{n})}(^{n}%
c_{0};F),
\]
where $q_{1}=...=q_{k}=2$ and $q_{k+1}=....=q_{n}=1.$ But, from{
\ }\cite[Theorem 8]{studia} we have%
\[
\mathcal{L}(^{n}c_{0};F)\neq\mathcal{L}_{as(q;q_{1},....,q_{n})}(^{n}%
c_{0};F),
\]
regardless of the $q<2$ and $q_{1}=...=q_{n}\geq1.$

Another relevant question is whether Corollary \ref{g} can be improved to
$p<q$, i.e.,

\begin{itemize}
\item If $\cot F=q<\infty$ and $E_{1}$,..., $E_{k}$ are infinite dimensional
$\mathcal{L}_{\infty}$-spaces, is there some $p<q$ for which, regardless of
the Banach spaces $E_{k+1},...,E_{n},$
\[
\mathcal{L}(E_{1},...,E_{k},...,E_{n};F)=\mathcal{L}_{as(p;q_{1},....,q_{n}%
)}(E_{1},...,E_{k},...,E_{n};F),
\]
where $q_{1}=...=q_{k}=2$ and $q_{k+1}=....=q_{n}=1$?
\end{itemize}

Again, applying \cite[Theorem 8]{studia} we obtain a negative answer to this question.

\section{Weak compactness and connections between the different classes}

It is well known that every absolutely $p$-summing operator is weakly compact
and completely continuous. So, a natural question is to ask whether their
multilinear generalizations still preserve these properties. In this section
we will obtain certain inclusions related to the different classes
investigated in this paper and we apply our results to face the aforementioned
question. In particular, we give an alternative direct answer for a question
posed by V. Dimant \cite{Dimant} and recently answered by Carando-Dimant
\cite{Carando}.

\begin{theorem}
\label{ss}(i) $[\Pi_{as(p)}](^{n}E;F)=\mathcal{L}_{as(p;p,\infty,...,\infty
)}(^{n}E;F)\cap...\cap\mathcal{L}_{as(p;\infty,...,\infty,p)}(^{n}E;F)$.

(ii) $\mathcal{L}_{d,p}(^{n}E;F)\subset\mathcal{L}_{si,p}(^{n}E;F)\subset
\lbrack\Pi_{as(p)}](^{n}E;F)$.

(iii) $\mathcal{L}_{si,p}(^{n}E;F)\subset\mathcal{L}_{d,np}(^{n}E;F).$

(iv) $\mathcal{L}_{si,p}(^{n}E;F)\subset\mathcal{L}_{fas,p}(^{n}%
E;F)\subset\mathcal{L}_{as,p}^{ev}(^{n}E;F)\subset\mathcal{L}_{as,p}(^{n}E;F)$.

(v) $\mathcal{L}_{si,p}(^{n}E;F)\subset\mathcal{L}_{sas,p}(^{n}E;F).$
\end{theorem}

Proof. (i) The case $n=3$ is illustrative. If $T\in\lbrack\Pi_{as(p)}%
](^{3}E;F)$ and $(x_{j})_{j=1}^{\infty}\in l_{p}^{w}(E),(y_{j})_{j=1}^{\infty
}\in l_{\infty}(E),(z_{j})_{j=1}^{\infty}\in l_{\infty}(E)$ are non
identically null, we have
\begin{align*}
&  (\sum_{j=1}^{\infty}\left\Vert T(x_{j},y_{j},z_{j})\right\Vert ^{p}%
)^{\frac{1}{p}}\\
&  =\left\Vert (y_{j})_{j=1}^{\infty}\right\Vert _{\infty}\left\Vert
(z_{j})_{j=1}^{\infty}\right\Vert _{\infty}\left(  \sum_{j=1}^{\infty
}\left\Vert T(x_{j},\frac{y_{j}}{\left\Vert (y_{j})_{j=1}^{\infty}\right\Vert
_{\infty}},\frac{z_{j}}{\left\Vert (z_{j})_{j=1}^{\infty}\right\Vert _{\infty
}})\right\Vert ^{p}\right)  ^{\frac{1}{p}}\\
&  \leq\left\Vert (y_{j})_{j=1}^{\infty}\right\Vert _{\infty}\left\Vert
(z_{j})_{j=1}^{\infty}\right\Vert _{\infty}\left(  \sum_{j=1}^{\infty
}\left\Vert \Psi_{1}^{(3)}(T)(x_{j})\right\Vert ^{p}\right)  ^{\frac{1}{p}}\\
&  \leq\left\Vert \Psi_{1}^{(3)}(T)\right\Vert _{as,p}\left\Vert (y_{j}%
)_{j=1}^{\infty}\right\Vert _{\infty}\left\Vert (z_{j})_{j=1}^{\infty
}\right\Vert _{\infty}\left\Vert (x_{j})_{j=1}^{\infty}\right\Vert _{w,p},
\end{align*}
and thus $T\in\mathcal{L}_{as(p;p,\infty,\infty)}(^{3}E;F)$. The other cases
are similar. The converse is not difficult.

(ii) The proof of $\mathcal{L}_{si,p}(^{n}E;F)\subset\lbrack\Pi_{as(p)}%
](^{n}E;F)$ is a direct consequence of (iv) of Proposition \ref{demmm1}. If
$T\in\mathcal{L}_{d,p}(^{n}E;F)$, then
\begin{align*}
&  \parallel T(x_{1,j},...,x_{n,j})\parallel\leq C(\int_{B_{E_{1}^{\prime}}%
}\mid\varphi_{1}(x_{1,j})\mid^{p}d\mu_{1}(\varphi_{1}))^{1/p}...(\int
_{B_{E_{n}^{\prime}}}\mid\varphi_{n}(x_{n,j})\mid^{p}d\mu_{n}(\varphi
_{n}))^{1/p}\\
&  =C(\int_{B_{E_{1}^{\prime}}\times...\times B_{E_{n}^{\prime}}}\mid
\varphi_{1}(x_{1,j})...\varphi_{n}(x_{n,j})\mid^{p}d(\mu_{1}\otimes
...\otimes\mu_{n})(\varphi_{1},...\varphi_{n}))^{1/p}%
\end{align*}
and hence $T\in\mathcal{L}_{si,p}(^{n}E;F)$.

(iii) Suppose that $T$ is $p$-semi integral. Then, by Theorem \ref{02},
\begin{align*}
\left(  \sum\limits_{j=1}^{m}\parallel T(x_{1,j},...,x_{n,j})\parallel
^{p}\right)  ^{1/p}  &  \leq C\left(  \underset{l=1,...,n}{\underset
{\varphi_{l}\in B_{E_{l}^{\prime}}}{\sup}}\sum\limits_{j=1}^{m}\mid\varphi
_{1}(x_{1,j})...\varphi_{n}(x_{n,j})\mid^{p}\right)  ^{1/p}\\
&  \leq C\underset{l=1,...,n}{\underset{\varphi_{l}\in B_{E_{l}^{\prime}}%
}{\sup}}\left(  \sum\limits_{j=1}^{m}\mid\varphi_{1}(x_{1,j})\mid^{np}\right)
^{\frac{1}{np}}...\left(  \sum\limits_{j=1}^{m}\mid\varphi_{n}(x_{n,j}%
)\mid^{np}\right)  ^{\frac{1}{np}}\\
&  =C\left\Vert (x_{1,j})_{j=1}^{m}\right\Vert _{w,np}...\left\Vert
(x_{n,j})_{j=1}^{m}\right\Vert _{w,np},
\end{align*}
and we are done.

(iv) If $T\in\mathcal{L}_{si,p}(^{n}E;F),$ then
\[
\underset{j_{1},...j_{n}=1}{\overset{m}{\sum}}\parallel T(x_{1,j_{1}%
},...,x_{n,j_{n}})\parallel^{p}\leq C^{p}\int_{B_{E_{1}^{\prime}}%
\times...\times B_{E_{n}^{\prime}}}\underset{j_{1},...j_{n}=1}{\overset
{m}{\sum}}\mid\varphi_{1}(x_{1,j_{1}})...\varphi_{n}(x_{n,j_{n}})\mid^{p}%
d\mu(\varphi_{1},...,\varphi_{n}).
\]
So
\begin{align*}
&  \left(  \underset{j_{1},...j_{n}=1}{\overset{m}{\sum}}\parallel
T(x_{1,j_{1}},...,x_{n,j_{n}})\parallel^{p}\right)  ^{1/p}\\
&  \leq C\left(  \int_{B_{E_{1}}^{\prime}\times...\times B_{E_{n}}^{\prime}%
}\underset{j_{1},...j_{n}=1}{\overset{m}{\sum}}\mid\varphi_{1}(x_{1,j_{1}%
})...\varphi_{n}(x_{n,j_{n}})\mid^{p}d\mu(\varphi_{1},...,\varphi_{n})\right)
^{\frac{1}{p}}\\
&  \leq C\underset{l=1,...n}{\sup_{\varphi_{l}\in B_{E_{l^{\prime}}}}}\left(
\underset{j_{1},...j_{n}=1}{\overset{m}{\sum}}\mid\varphi_{1}(x_{1,j_{1}%
})...\varphi_{n}(x_{n,j_{n}})\mid^{p}\right)  ^{1/p}\\
&  =C\underset{l=1,...n}{\underset{\varphi_{l}\in B_{E_{l^{\prime}}}}{\sup}%
}\left[  \underset{j_{1}=1}{(\overset{m}{\sum}}\mid\varphi_{1}(x_{1,j_{1}%
})\mid^{p})^{1/p}...\underset{j_{n}=1}{(\overset{m}{\sum}}\mid\varphi
_{n}(x_{n,j_{n}})\mid^{p})^{1/p}\right]  ,
\end{align*}
and thus $T\in\mathcal{L}_{fas,p}(^{n}E;F).$ Now let us consider
$T\in\mathcal{L}_{fas,p}(^{n}E;F)$. The case $n=2$ is illustrative and
indicates the proof. If $(x_{j})_{j=1}^{\infty},(y_{j})_{j=1}^{\infty}\in
l_{p}^{w}(E)$, we have%
\begin{align*}
&  \left(  \sum\limits_{j=1}^{\infty}\parallel T(a+x_{j},b+y_{j}%
)-T(a,b)\parallel^{p}\right)  ^{\frac{1}{p}}\\
&  \leq\left(  \sum\limits_{j=1}^{\infty}\parallel T(a,y_{j})\parallel
^{p}\right)  ^{\frac{1}{p}}+\left(  \sum\limits_{j=1}^{\infty}\parallel
T(x_{j},b)\parallel^{p}\right)  ^{\frac{1}{p}}+\left(  \sum\limits_{j=1}%
^{\infty}\parallel T(x_{j},y_{j})\parallel^{p}\right)  ^{\frac{1}{p}}\\
&  \leq\left(  \sum\limits_{j,k=1}^{\infty}\parallel T(z_{k},y_{j}%
)\parallel^{p}\right)  ^{\frac{1}{p}}+\left(  \sum\limits_{j,k=1}^{\infty
}\parallel T(x_{j},w_{k})\parallel^{p}\right)  ^{\frac{1}{p}}+\left(
\sum\limits_{j,k=1}^{\infty}\parallel T(x_{j},y_{k})\parallel^{p}\right)
^{\frac{1}{p}}<\infty,
\end{align*}
where $(z_{j})_{j=1}^{\infty}=(a,0,0,...)$ and $(w_{j})_{j=1}^{\infty
}=(b,0,0,...).$

(v) If $T$ is $p$-semi integral, then
\begin{align*}
\left(  \underset{j=1}{(\overset{m}{\sum}}\parallel T(x_{1,j},...,x_{n,j}%
)\parallel^{p}\right)  ^{1/p}  &  \leq C\left(  \underset{_{\varphi_{l}\in
B_{E_{l^{\prime}}}}}{\sup}\sum\limits_{j=1}^{m}\mid\varphi_{1}(x_{1,j}%
)...\varphi_{n}(x_{n,j})\mid^{p}\right)  ^{1/p}\\
&  \leq C\left(  \underset{\phi\in B_{\mathcal{L}(E_{1},...,E_{n};\mathbb{K}%
)}}{\sup}\sum\limits_{j=1}^{m}\mid\phi(x_{1,j},...,x_{n,j})\mid^{p}\right)
^{1/p}%
\end{align*}
and thus $T\in\mathcal{L}_{sas,p}(^{n}E;F)$.

\begin{remark}
Obviously, each one of the assertions of the Theorem \ref{ss} holds for spaces
$E_{1},...,E_{n}$ instead of $E,...,E.$ The inclusion $\mathcal{L}%
_{si,1}(E_{1},...E_{n};F)\subset\mathcal{L}_{sas,1}(E_{1},...E_{n};F)$ is
strict. In fact, if $T:l_{2}\times l_{2}\rightarrow\mathbb{K}$ is given by
\[
T((x)_{j=1}^{\infty},(y_{j})_{j=1}^{\infty})=\sum\limits_{j=1}^{\infty}%
y_{j}\sum\limits_{k=1}^{j}x_{k},
\]
then $T$ fails to be semi integral (see \cite{AlencarMatos})$,$ but $T$ is
strongly $1$-summing, because every continuous $n$-linear form is obviously
strongly $1$-summing.

The inclusion $\mathcal{L}_{si,1}(E_{1},...E_{n};F)\subset\mathcal{L}%
_{fas,1}(E_{1},...E_{n};F)$ is also strict, since
\[
\mathcal{L}(l_{2},l_{1};\mathbb{K})=\mathcal{L}_{fas,1}(l_{2},l_{1}%
;\mathbb{K})\text{ (\cite{DM})}%
\]
and
\[
\mathcal{L}(l_{2},l_{1};\mathbb{K})\neq\mathcal{L}_{si,1}(l_{2},l_{1}%
;\mathbb{K}).
\]
In fact, if we had $\mathcal{L}(l_{2},l_{1};\mathbb{K})=\mathcal{L}%
_{si,1}(l_{2},l_{1};\mathbb{K}),$ we would obtain
\[
\mathcal{L}(l_{2},l_{\infty})=\mathcal{L}_{as,1}(l_{2},l_{\infty})
\]
and it is a contradiction. The inclusion
\[
\mathcal{L}_{fas,1}(E_{1},...E_{n};F)\subset\mathcal{L}_{as,1}(E_{1}%
,...E_{n};F)\text{ }%
\]
is also strict (see \cite{Collect}).

It is interesting to observe that (in general) $\mathcal{L}_{sas,p}%
(E_{1},...E_{n};F)$ is not contained in $\mathcal{L}_{as,p}(E_{1},...E_{n};F)$
and $\mathcal{L}_{as,p}(E_{1},...E_{n};F)$ is not contained in $\mathcal{L}%
_{sas,p}(E_{1},...E_{n};F).$ In fact, $\mathcal{L}_{as,1}(^{2}l_{1}%
;l_{1})=\mathcal{L}(l_{1};l_{1})$ and $\mathcal{L}_{sas,1}(^{2}l_{1}%
;l_{1})\neq\mathcal{L}(^{2}l_{1};l_{1}).$ On the other hand $\mathcal{L}%
_{sas,2}(^{2}l_{2};\mathbb{K})=\mathcal{L}(^{2}l_{2};\mathbb{K})$ and
$\mathcal{L}_{as,2}(^{2}l_{2};\mathbb{K})\neq\mathcal{L}(^{2}l_{2}%
;\mathbb{K}).$
\end{remark}

The next result shows that the spaces of semi integral and dominated mappings
coincides in some situations:

\bigskip

\begin{theorem}
\label{aa}(i) If $\cot E=2,$ then $\mathcal{L}_{si,1}(^{2}E;F)=\mathcal{L}%
_{d,1}(^{2}E;F)$ for every Banach space $F$.

(ii) If $E$ is an $\mathcal{L}_{\infty}$-space, then $\mathcal{L}_{si,1}%
(^{n}E;F)=\mathcal{L}_{d,1}(^{n}E;F)$ for every $n$ and every $F$.
\end{theorem}

Proof. (i) If $E$ has cotype $2$, we know that $\mathcal{L}_{as,1}%
(E;F)=\mathcal{L}_{as,2}(E;F)$ for every Banach space $F$. Thus, if
$T\in\mathcal{L}_{d,2}(^{2}E;F),$ then $T=R(u_{1},u_{2}),$ with $R\in
\mathcal{L}(^{2}G;F)$ and $u_{1},u_{2}\in\mathcal{L}_{as,2}(E;G)=\mathcal{L}%
_{as,1}(E;G).$ Hence $T\in\mathcal{L}_{d,1}(^{2}E;F)$ and thus
\[
\mathcal{L}_{d,1}(^{2}E;F)=\mathcal{L}_{d,2}(^{2}E;F)
\]
for every $F$. Since $\mathcal{L}_{d,1}(^{2}E;F)\subset\mathcal{L}_{si,1}%
(^{2}E;F)\subset\mathcal{L}_{d,2}(^{2}E;F)$ we thus have $\mathcal{L}%
_{si,1}(^{2}E;F)=\mathcal{L}_{d,1}(^{2}E;F).$

For the proof of (ii), a localization argument allows to consider $E=C(K),$
where $K$ is a compact Hausdorff space. By applying \cite[Proposition
2.6]{Villa}, it is not hard to see that every Pietsch-integral $n$-linear
mapping is $1$-dominated. Besides, the Theorem \ref{ss} asserts that every
$1$-dominated mapping is $1$-semi-integal. On the other hand, every
$1$-semi-integral mapping on $C(K)$ is Pietsch integral \cite[Theorem
5.6]{AlencarMatos}and the proof is done.

\begin{remark}
It is known that $\mathcal{L}_{d,p}(^{n}E;F)\neq\lbrack\Pi_{as(p)}](^{n}E;F).$
In fact, the $2$-linear mapping $T:l_{2}\times l_{2}\rightarrow\mathbb{K}%
\ $given by $T(x,y)=\sum\limits_{j=1}^{\infty}\frac{1}{j^{\alpha}}x_{j}y_{j}$
with $\alpha=\frac{1}{2}+\varepsilon$ and $0<\varepsilon<\frac{1}{2}$ is so
that $T\in$ $[\Pi_{as(1)}](^{2}l_{2};\mathbb{K})$ and fails to be
$1$-dominated (see \cite{Pellegrino})$.$ This example and the last proposition
shows that in general $\mathcal{L}_{si,1}(^{n}E;F)\neq\lbrack\Pi_{as(1)}%
](^{n}E;F).$
\end{remark}

In \cite{Botelho2}, Botelho proved that $P_{n}:l_{1}\rightarrow l_{1}%
:P_{n}((\alpha_{i})_{i=1}^{\infty})=((\alpha_{i})^{n})_{i=1}^{\infty}$ was
$n$-dominated and was not weakly compact. The same occurs with the symmetric
$n$-linear mapping associated to $P$.

The question \textquotedblleft Is every strongly $p$-summing $n$-linear
mapping weakly compact?\textquotedblright\ appears in \cite{Dimant} and was
recently answered by Carando-Dimant in \cite{Carando}. However, by Theorem
\ref{ss}, since $\mathcal{L}_{d,p}(E_{1},...E_{n};F)\subset\mathcal{L}%
_{sas,p}(E_{1},...E_{n};F),$ one can realize that Botelho$^{\text{'}}$s
counterexample is a (more general) answer to this question.

Concerning completely continuous mappings, it is not hard to prove that every
continuous $p$-semi integral mapping is completely continuous. On the other
hand, contrary to the linear case, the absolutely summing (and strongly
summing) multilinear mappings are not completely continuous, in general. For
example, $T:l_{2}\times l_{2}\rightarrow\mathbb{K}$ given by\ $T((x_{n}%
)_{n=1}^{\infty},(y_{n})_{n=1}^{\infty})=\sum\limits_{j=1}^{\infty}x_{j}y_{j}%
$\ is absolutely $1$-summing and strongly $2$-summing but fails to be
completely continuous.\

\section{A new class}

In this section, we introduce a new class related to the concept of absolute
summability. Our idea is to join two interesting ideas: to sum in multiple
index and to work with multilinear mappings instead of linear functionals, as
in the definitions of fully (multiple) summing mappings and strongly summing
mappings, respectively.

We will say that $T\in\mathcal{L}(E_{1},...E_{n};F)$ is strongly fully
$p$-summing if and only if there exists $C\geq0$ such that
\begin{equation}
\left(  \sum\limits_{j_{1},...,j_{n}=1}^{m}\parallel T(x_{j_{1}}%
^{(1)},...,x_{j_{n}}^{(n)})\parallel^{p}\right)  ^{1/p}\leq C\left(
\underset{}{\underset{\phi\in B_{\mathcal{L}(E_{1},...,E_{n})}}{\sup}}%
\sum\limits_{j_{1},...,j_{n}=1}^{m}\mid\phi(x_{j_{1}}^{(1)},...,x_{j_{n}%
}^{(n)})\mid^{p}\right)  ^{1/p}%
\end{equation}
for every natural $m$, $x_{j_{l}}^{(l)}\in E_{l}$ with $l=1,...,n$ and
$j=1,...,m.$ The space of all strongly fully $p$-summing $n$-linear mappings
from $E_{1}\times...\times E_{n}$ into $F$ will be denoted by $\mathcal{L}%
_{sf,p}(E_{1},...,E_{n};F)$ and the infimum of the $C$ for which the
inequality always holds defines a norm $\Vert.\Vert_{sf,p}$ for $\mathcal{L}%
_{sf,p}(E_{1},...,E_{n};F)$. Under this norm, $\mathcal{L}_{sf,p}%
(E_{1},...,E_{n};F)$ is complete. One can verify the following properties:

(i) $\mathcal{L}_{fas,p}(E_{1},...,E_{n};F)\subset\mathcal{L}_{sf,p}%
(E_{1},...,E_{n};F).$

(ii) $\mathcal{L}_{sas,p}(E_{1},...,E_{n};F)\subset\mathcal{L}_{sf,p}%
(E_{1},...,E_{n};F).$

(iii) If $\mathcal{L}_{sf,p}(E_{1},...,E_{n};F)=\mathcal{L}(E_{1}%
,...,E_{n};F),$ then $\mathcal{L}(E_{j_{1}},...E_{j_{n}};F)=\mathcal{L}%
_{sf,p}(E_{j_{1}},...E_{j_{n}};F)$, for every $j_{1},...,j_{k}$ in
$\{1,...,n\}$ with $j_{r}\neq j_{s}$ if $r\neq s.$

(iv) $\mathcal{L}_{sf,p}(^{n}E;E)=\mathcal{L}(^{n}E;E)\Longleftrightarrow\dim
E<\infty.$

(v) If $T\in\mathcal{L}(E_{1},...,E_{n};F),$ then the Aron-Berner extension of
$T$ belongs to $\mathcal{L}(E_{1}^{\prime\prime},...,E_{n}^{\prime\prime};F)$.

Since $\mathcal{L}_{sf,p}(E_{1},...,E_{n};F)$ contains $\mathcal{L}%
_{fas,p}(E_{1},...,E_{n};F)$ and $\mathcal{L}_{sas,p}(E_{1},...,E_{n};F)$,
every coincidence result for strongly $p$-summing and/or fully $p$-summing
multilinear mappings still holds for the strongly fully summing mappings. On
the other hand, (iv) implies that $\mathcal{L}_{sf,p}(E_{1},...,E_{n};F)$ has
a Dvoretzky-Rogers type theorem and (iii) shows that coincidence results
$\mathcal{L}_{sf,p}(E_{1},...,E_{n};F)=\mathcal{L}(E_{1},...,E_{n};F)$ are not
so common. For example, since $\mathcal{L}(^{n}l_{1};l_{2})=\mathcal{L}%
_{sas,1}(^{n}l_{1};l_{2})$ (\cite{Dimant}), by (ii), (iii)  and
\cite[Theorem 4.2]{Lindenstrauss} we can prove that if $E$ has
unconditional Schauder basis and $F$ is an infinite dimensional
Banach space, then
$\mathcal{L}_{sf,p}(^{n}E;F)=\mathcal{L}(^{n}E;F)$ if and only if
$E$ is isomorphic to $l_{1}$ and $F$ is a Hilbert space.



\begin{thebibliography}{99}                                                                                               %


\bibitem {AlencarMatos}R. Alencar and M.C. Matos, Some classes of multilinear
mappings between Banach spaces, Publicaciones Departamento An\'{a}lisis
Matematico, Universidad Complutense Madrid, Section 1, number 12 (1989).

\bibitem {Perez}F. Bombal, D. P\'{e}rez-Garc\'{\i}a and I. Villanueva,
Multilinear extensions of a Grothendieck's theorem, preprint.

\bibitem {Botelhoirish}G. Botelho, Cotype and absolutely summing multilinear
mappings and homogeneous polynomials, Proc. Roy. Irish Acad. Sect. A,
\textbf{97} (1997), 145-153.

\bibitem {Botelho2}G. Botelho, Weakly compact and absolutely summing
polynomials, J. Math. Anal. Appl. \textbf{265} (2002), 458-462.

\bibitem {BP}G. Botelho and D.M. Pellegrino, Dominated polynomials in
$\mathcal{L}_{p}$-spaces, preprint.

\bibitem {Carando}D. Carando and V. Dimant, On summability of bilinear
operators, Math. Nachr. \textbf{259} (2003), 3-11.

\bibitem {Diestel}J. Diestel, H. Jarchow and A. Tonge, \textit{Absolutely
Summing Operators}, Cambridge Stud.\ Adv.\ Math. \textbf{43}, Cambridge
University Press, Cambridge 1995.

\bibitem {Dimant}V. Dimant, Strongly $p$-summing multilinear mappings, J.
Math. Anal. Appl. \textbf{278} (2003), 182-193.

\bibitem {Lindenstrauss}J. Lindenstrauss and A. Pe\l czy\'{n}ski, Absolutely
summing operators in $\mathcal{L}_{p}$-spaces and their applications, Studia
Math. \textbf{29} 1968, 275-324.

\bibitem {Anais}M.C. Matos, Absolutely summing holomorphic mappings, An. Acad.
bras. Ci., \textbf{68} (1996), 1-13.

\bibitem {Collect}M.C. Matos, Fully absolutely summing mappings and Hilbert
Schmidt operators, Collect. Mat. \textbf{54} (2003), 111-136.

\bibitem {Nach}M.C. Matos, Nonlinear absolutely summing mappings, Math. Nachr.
\textbf{258} (2003), 71-89.

\bibitem {Melendez}Y. Mel\'{e}ndez and A. Tonge, Polynomials and the Pietsch
domination Theorem, Proc. Roy. Irish Acad. Sect. A (1999), 195-212.

\bibitem {studia}D.M. Pellegrino, Cotype and absolutely summing homogeneous
polynomials in $\mathcal{L}_{p}$ spaces, Studia Math. \textbf{157}, (2003), 121-131.

\bibitem {Pellegrino}D.M. Pellegrino, Aplica\c{c}\~{o}es entre espa\c{c}os de
Banach relacionadas \`{a} converg\^{e}ncia de s\'{e}ries, Thesis, Unicamp 2002.

\bibitem {cotype}D.M. Pellegrino, Cotype and nonlinear absolutely summing
mappings, preprint.

\bibitem {DM}D.M. Pellgrino and M.L.V. Souza, Fully summing multilinear and
holomorphic mappings into HIlbert spaces, to appear in Math. Nachr.

\bibitem {PG}D. P\'{e}rez-Garc\'{\i}a, Operadores Multilineales absolutamente
sumantes, Dissertation, Universidad Complutense de Madrid 2002.

\bibitem {P1967}A. Pietsch, Absolute $p$-summierende Abbildungen in normierten
R\"{a}umen, Studia Math. \textbf{28} (1967), 333-353.

\bibitem {Schneider}B. Schneider, On absolutely $p$-summing and related
multilinear mappings. Wissenchaftliche Zeitschrift der Brandenburger
Landeshochschule \textbf{35} (1991), 105-117.

\bibitem {Souza}M.L.V. Souza, Aplica\c{c}\~{o}es multilineares completamente
absolutamente somantes, in ``55 Semin\'{a}rio Brasileiro de An\'{a}lise'', May
2002, 689-706.

\bibitem {CA}C.A. Soares, Aplica\c{c}\~{o}es multilineares e polin\^{o}mios
misto somantes, Thesis, Unicamp 1998.

\bibitem {Villa}I. Villanueva, Integral mappings between Banach spaces, J.
Math. Anal. Appl. \textbf{279} (2003), 56-70.


\end{thebibliography}
\end{document}